\documentclass{article}
\usepackage{epsf,epsfig}

\newcommand{\Chi}{X}

\begin{document}

\title{Eigenmodes of Lens and Prism Spaces}

\author{Roland Lehoucq \\
    {\it\small CE-Saclay, DSM/DAPNIA/Service d'Astrophysique} \\
    {\it\small F-91191 Gif sur Yvette Cedex (France)} \\
    {\it\small roller@discovery.saclay.cea.fr} \\
    \\
    Jean-Philippe Uzan \\
    {\it\small Institut d'Astrophysique de Paris} \\
    {\it\small 98 bis, Bd Arago, 75014 Paris (France)} \\
    {\it\small Laboratoire de Physique Th\'eorique, CNRS-UMR 8627, B\^at. 210} \\
    {\it\small Universit\'e Paris XI, F--91405 Orsay Cedex (France)} \\
    {\it\small uzan@amorgos.unige.ch} \\
    \\
    Jeffrey Weeks \\
    {\it\small 15 Farmer Street} \\
    {\it\small Canton NY 13617-1120 (USA)} \\
    {\it\small weeks@northnet.org}}

\maketitle

\begin{abstract}
    Cosmologists are taking a renewed interest in multiconnected
    spherical 3-mani\-folds (spherical spaceforms) as possible models
    for the physical universe.
    To understand the formation of large scale structures in such
    a universe, cosmologists express physical quantities, such as density
    fluctuations in the primordial plasma, as linear combinations
    of the eigenmodes of the Laplacian, which can then be integrated
    forward in time.
    This need for explicit eigenmodes contrasts sharply with
    previous mathematical investigations, which have focused
    on questions of isospectrality rather than eigenmodes.
    The present article provides explicit orthonormal bases
    for the eigenmodes of lens and prism spaces.
\end{abstract}

\pagebreak

\section{Introduction}

In recent years cosmologists have taken a renewed interest
in multiply connected 3-manifolds as possible models for
the universe \cite{LuminetStarkmanWeeks,CornishWeeks,Luminet},
motivated by upcoming opportunities to determine the topology
of the real universe using satellite measurements
of the microwave background \cite{CornishSpergelStarkman}
and galaxy catalogs \cite{LLL}.
Cosmologists initially focused on closed hyperbolic 3-manifolds,
favored by the low observed matter density in the universe,
as well as the more easily understood flat 3-manifolds.
But since 1998 it has become clear that the modest amount
of matter in the universe (30\%) is complemented by a large
amount of exotic energy (70\%).
This extra energy implies that the observable universe
is approximately flat, or perhaps slightly spherical \cite{JaffeEtAl}.
Cosmologists have therefore shifted their interest from
hyperbolic 3-manifolds to flat and spherical ones.
Beyond the data's very slight preference for a spherical universe,
the cosmologists' new interest in multiply connected spherical
3-manifolds (spherical spaceforms) is due to the fact that
the volume of a spherical 3-manifold decreases as the
topology gets more complicated, unlike hyperbolic 3-manifolds
whose volumes increase as the topology gets more complicated.
Thus even though both hyperbolic and spherical 3-manifolds
are consistent with an approximately flat observable universe,
the spherical topologies would be more easily
detectable observationally \cite{GLLUW}.

To understand and simulate microwave background measurements
in a multiply connected spherical universe, cosmologists must
first understand the density fluctuations in the primordial plasma
(see \cite{Levin} for a review).
Such density fluctuations are expressed as linear combinations
of the eigenmodes (eigenfunctions) of the Laplacian, just as
the vibration of a drumhead may be expressed as a linear combination
of the drumhead's eigenmodes.
But unlike mathematicians' studies of isospectrality \cite{IkedaYamamoto},
where the spectrum was the primary object of interest
(``Can you hear the shape of a lens space?'')
and the eigenmodes were secondary, the cosmologists'
research puts the eigenmodes at center stage.  More specifically,
for each wave number $k$ (corresponding to eigenvalue $k(k+2)$),
the cosmologists want an explicit orthonormal basis
for the corresponding space of eigenmodes.
The present paper provides such an eigenbasis for all
lens spaces (Theorem 2, Section \ref{SectionLensSpaces})
and prism spaces (Theorem 3, Section \ref{SectionPrismSpaces}).

\section{Toroidal coordinates}

\begin{figure}
\begin{center}
\epsfig{file=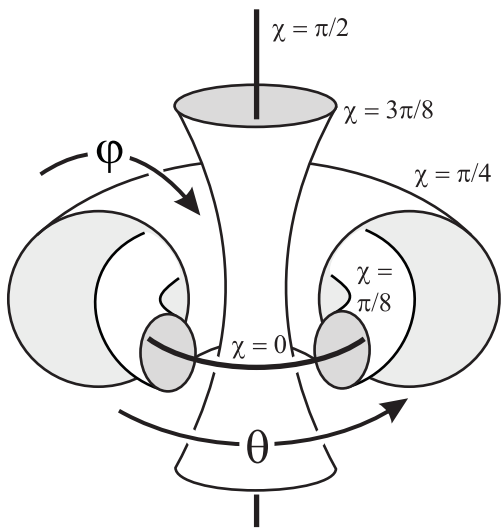}
\caption{Toroidal coordinates.
Nested tori fill the 3-sphere like layers of an onion.
Just as the layers of an onion collapse to a line at the onion's core,
the nested tori collapse to two circles, one at $\chi = 0$
and the other at $\chi = \pi/2$.} \label{fig1}
\end{center}
\end{figure}

The determination of the eigenmodes of a lens space (resp. prism space)
is elementary and constructive.  Visualize such a manifold as
the 3-sphere $S^3$ under the action of a cyclic (resp. binary dihedral)
group $\Gamma$ of covering transformations.  The key to simplicity
is to choose a coordinate system that respects the covering
transformations $\Gamma$.  A {\it toroidal coordinate system}
meets our needs perfectly (Figure~\ref{fig1}).

Let $x$, $y$, $z$, and $w$ be the usual coordinates in ${\bf R}^4$,
so the 3-sphere $S^3$ is defined by $x^2 + y^2 + z^2 + w^2 = 1$.
The coordinates $\chi$, $\theta$, and $\varphi$ parameterize the 3-sphere
as
\begin{equation}
\begin{array}{ccl}
x &=& \cos{\chi}\,\cos{\theta} \\
y &=& \cos{\chi}\,\sin{\theta} \\
z &=& \sin{\chi}\,\cos{\varphi}\\
w &=& \sin{\chi}\,\sin{\varphi}
\label{CoordinateDefinition}
\end{array}
\end{equation}
for
\begin{equation}
\begin{array}{ccccl}
0 & \leq &\chi&    \leq & \pi/2 \\
0 & \leq &\theta&  \leq & 2\pi  \\
0 & \leq &\varphi& \leq & 2\pi.
\end{array}
\end{equation}
For each fixed value of $\chi \in (0, \pi/2)$, the $\theta$ and
$\varphi$ coordinates sweep out a torus.  Taken together,
these tori almost fill $S^3$.  The exceptions occur
at the endpoints $\chi = 0$ and $\chi = \pi/2$, where the
stack of tori collapses to the circles $x^2 + y^2 = 1$ and
$z^2 + w^2 = 1$, respectively.

\section{The Laplacian in toroidal coordinates}

The coordinates $\chi$, $\theta$, and $\varphi$ are everywhere
orthogonal to each other.  Thus the metric on the 3-sphere
may be written as
\begin{equation}
ds^2 = h_{\chi}^2\,d\chi^2 + h_{\theta}^2\,d\theta^2 + h_{\varphi}^2\,d\varphi^2
\end{equation}
where
\begin{equation}
\begin{array}{ccl}
\label{hvalues}
h_{\chi}    &=& 1           \\
h_{\theta}  &=& \cos{\chi}  \\
h_{\varphi} &=& \sin{\chi}.
\end{array}
\end{equation}
The Laplacian is just the divergence of the gradient of a function $\Psi$,
and the divergence is, in turn, just the ``net outflow per unit volume''.
So by visualizing the small volume element $d\chi\,d\theta\,d\varphi$
we may write down the Laplacian for any orthogonal coordinate system as
\begin{equation}\label{GeneralLaplacian}
\nabla^2 = \frac{1}{h_{\chi}h_{\theta}h_{\varphi}}
    \left\{
        \frac{\partial}{\partial \chi}
            \frac{h_{\theta}h_{\varphi}}{h_{\chi}}
            \frac{\partial}{\partial \chi}
       +
         \frac{\partial}{\partial \theta}
            \frac{h_{\chi}h_{\varphi}}{h_{\theta}}
            \frac{\partial}{\partial \theta}
       +
         \frac{\partial}{\partial \varphi}
            \frac{h_{\chi}h_{\theta}}{h_{\varphi}}
            \frac{\partial}{\partial \varphi}
  \right\}
\end{equation}
with no calculation required.  In the present case, substituting
(\ref{hvalues}) into (\ref{GeneralLaplacian}) gives the Laplacian
in toroidal coordinates
\begin{equation}\label{TorusLaplacian}
\nabla^2 = \frac{1}{\cos{\chi}\sin{\chi}}
    \left\{
        \frac{\partial}{\partial \chi}
            \cos{\chi}\sin{\chi}
            \frac{\partial}{\partial \chi}
       +
         \frac{\partial}{\partial \theta}
            \frac{\sin{\chi}}{\cos{\chi}}
            \frac{\partial}{\partial \theta}
       +
         \frac{\partial}{\partial \varphi}
            \frac{\cos{\chi}}{\sin{\chi}}
            \frac{\partial}{\partial \varphi}
   \right\}.
\end{equation}

\section{The Helmholtz equation in toroidal coordinates}

The wave number $k$ parameterizes the eigenmodes of the Laplacian
on the 3-sphere $S^3$.  Each integer wave number $k > 0$ corresponds
to an eigenvalue $-k(k+2)$ with multiplicity $(k+1)^2$
\cite{BergerGauduchonMazet,SteinWeiss}.
The Helmholtz equation thus takes the form
\begin{equation}\label{HelmholtzEquation}
\nabla^2\Psi = -k(k+2) \Psi.
\end{equation}
We will look for solutions that factor as
\begin{equation}\label{factorization}
\Psi(\chi,\theta,\varphi) = \Chi(\chi)\,\Theta(\theta)\,\Phi(\varphi).
\end{equation}
We have no {\it a priori} guarantee that all solutions must
take this form, but in Section \ref{BasisSection} we'll see that the number
of independent solutions of this form does indeed equal the
dimension $(k+1)^2$ of the full eigenspace.

Substituting the expression (\ref{TorusLaplacian}) for $\nabla^2$
and the factorization (\ref{factorization}) of $\Psi$ into the
Helmholtz equation (\ref{HelmholtzEquation}) gives
\begin{equation}
     \frac{\Theta\Phi}{\cos{\chi}\sin{\chi}}
     \frac{\partial}{\partial \chi}
     \cos{\chi}\sin{\chi}
     \frac{\partial \Chi}{\partial \chi}
  +
     \frac{\Chi\Phi}{\cos^2{\chi}}
     \frac{\partial^2\Theta}{\partial \theta^2}
  +
     \frac{\Chi\Theta}{\sin^2{\chi}}
     \frac{\partial^2\Phi}{\partial \varphi^2}
  =
     -k(k+2)\,\Chi\Theta\Phi.
\end{equation}
Multiplying through by $\cos^2{\chi}\sin^2{\chi}/(\Chi\Theta\Phi)$
isolates the $\Theta$ and $\Phi$ factors
\begin{eqnarray}
    \frac{\cos{\chi}\sin{\chi}}{\Chi}
     \frac{d}{d \chi}
     \cos{\chi}\sin{\chi}
     \frac{d \Chi}{d \chi}
  +
     \sin^2{\chi}
     (\frac{1}{\Theta}
     \frac{d^2\Theta}{d \theta^2})
  +
     \cos^2{\chi}
     (\frac{1}{\Phi}
     \frac{d^2\Phi}{d \varphi^2})  \nonumber \\
 =
     -k(k+2)\,\cos^2{\chi}\sin^2{\chi}.
\label{SeparatedVariables}
\end{eqnarray}
The expressions in $\Theta$ and $\Phi$ must each be constant,
and to allow a periodic solution the constants must be negative,
\begin{eqnarray}
\label{ThetaEquation}
\frac{1}{\Theta}\,\frac{d^2\Theta}{d \theta^2 } &=& -\ell^2 \\
\label{PhiEquation}
\frac{1}{\Phi  }\,\frac{d^2\Phi  }{d \varphi^2} &=& -m^2.
\end{eqnarray}
The solutions are the usual circular harmonics
\begin{equation}
\Theta_\ell(\theta) = \cos{|\ell| \theta }\quad{\rm or}\quad \sin{|\ell| \theta }
\label{ThetaSolution}
\end{equation}
and
\begin{equation}
\Phi_m(\varphi)   = \cos{|m| \varphi}\quad{\rm or}\quad \sin{|m| \varphi}.
\label{PhiSolution}
\end{equation}
By convention, nonnegative $\ell$ indicates $\cos{|\ell| \theta }$
while negative $\ell$ indicates $\sin{|\ell| \theta }$, and
similarly for $m$.

Substituting (\ref{ThetaEquation}) and (\ref{PhiEquation})
into the Helmholtz equation (\ref{SeparatedVariables}) reduces it
to a second order ordinary differential equation for $\Chi$
\begin{equation}
     \frac{\cos{\chi}\sin{\chi}}{\Chi}
     \frac{d}{d \chi}
     \cos{\chi}\sin{\chi}
     \frac{d \Chi}{d \chi}
  -
     \ell^2\,\sin^2{\chi}
  -
     m^2\,\cos^2{\chi}
  =
     -k(k+2)\,\cos^2{\chi}\sin^2{\chi}.
\label{RadialEquation}
\end{equation}
For integers $k$, $\ell$, and $m$ satisfying
$|\ell| + |m| \leq k$ and $\ell + m \equiv k$ (mod 2),
equation (\ref{RadialEquation}) is a close relative
of the Jacobi equation and admits the solution
\begin{equation}
\Chi_{k\ell m}(\chi) = \cos^{|\ell|}\chi\,\sin^{|m|}\chi\;
  P^{(|m|,|\ell|)}_d (\cos 2\chi)
\label{RadialSolution}
\end{equation}
where $P^{(|m|,|\ell|)}_d$ is the Jacobi polynomial
\begin{equation}
P^{(|m|,|\ell|)}_d(u) = \frac{1}{2^d}\;
    \sum_{i=0}^d {{|m|+d}\choose{i}} {{|\ell|+d}\choose{d-i}}
                   \,(u + 1)^i\,(u - 1)^{d-i}
\label{JacobiPolynomial}
\end{equation}
and
\begin{equation}
 d = \frac{k - (|\ell| + |m|)}{2}.
\end{equation}

\section{The eigenmodes of $S^3$ in trigonometric and polynomial forms}
\label{SectionTrigAndPolynomialForms}

Substituting the expressions for $\Chi$, $\Theta$, and $\Phi$
from (\ref{RadialSolution}), (\ref{ThetaSolution}), and
(\ref{PhiSolution}) gives the eigenmode
\begin{equation}
\begin{array}{rl}
\widetilde\Psi_{k\ell m}(\chi,\theta,\varphi) &
=
\cos^{|\ell|}\chi\,\sin^{|m|}\chi\;P^{(|m|,|\ell|)}_d (\cos 2\chi)  \\
 & \quad \times \quad (\cos{|\ell| \theta }\quad{\rm or}\quad \sin{|\ell| \theta })  \\
 & \quad \times \quad (\cos{|m| \varphi}\quad{\rm or}\quad \sin{|m| \varphi})
\end{array}
\label{PsiSolution}
\end{equation}
where as usual the choice of $\cos{|\ell| \theta }$ or $\sin{|\ell| \theta }$
(resp. $\cos{|m| \varphi}$ or $\sin{|m| \varphi}$) depends on
the sign of $\ell$ (resp. $m$).

The Jacobi polynomial may be expanded as a homogeneous polynomial
of degree $2d$ in $x$, $y$, $z$, and $w$,
\begin{eqnarray}
P^{(|m|,|\ell|)}_d(\cos 2\chi) &=& \frac{1}{2^d}\;
    \sum_{i=0}^d {{|m|+d}\choose{i}} {{|\ell|+d}\choose{d-i}}
                   \,(\cos 2\chi + 1)^i\,(\cos 2\chi - 1)^{d-i} \nonumber \\
&=& \sum_{i=0}^d {{|m|+d}\choose{i}} {{|\ell|+d}\choose{d-i}}
                   \,(\cos^2\chi)^i\,(-\sin^2 \chi)^{d-i} \nonumber \\
&=& \sum_{i=0}^d {{|m|+d}\choose{i}} {{|\ell|+d}\choose{d-i}}
                   \,(x^2 + y^2)^i\,(-(z^2 + w^2))^{d-i}.
\end{eqnarray}
The $\cos^{|\ell|}\chi$ factor combines felicitously with
the $\cos |\ell| \theta$ or $\sin |\ell| \theta$ factor to create
a homogeneous polynomial of degree $|\ell|$ in $x$ and $y$, for example,
\begin{eqnarray}
& & \cos^{|\ell|}\chi\;\cos |\ell| \theta \nonumber \\
&=& \cos^{|\ell|}\chi\; \sum_{0 \leq i \leq |\ell|/2}
  (-1)^i {|\ell| \choose 2i} \cos^{|\ell| - 2i}\theta \sin^{2i}\theta \nonumber \\
&=& \sum_{0 \leq i \leq |\ell|/2} (-1)^i {|\ell| \choose 2i}
  (\cos\chi\,\cos\theta)^{|\ell| - 2i} (\cos\chi\,\sin\theta)^{2i} \nonumber \\
&=& \sum_{0 \leq i \leq |\ell|/2} (-1)^i {|\ell| \choose 2i}
  x^{|\ell| - 2i} y^{2i}.
\end{eqnarray}
Similarly, $\sin^{|m|}\chi$ combines with
$\cos |m| \varphi$ or $\sin |m| \varphi$ to create
a degree $|m|$ polynomial in $z$ and $w$.
Multiplied together, these factors express
$\widetilde\Psi_{k\ell m}$ as a homogeneous degree $k$
harmonic polynomial in $(x,y,z,w)$ coordinates.
For example, when $k = 7$, $\ell = 3$, and $m = -2$ we have
\begin{eqnarray}
\widetilde\Psi_{7,3,-2}(\chi,\theta,\varphi)
&=& [ P^{(2,3)}_1 (\cos 2\chi) ] \;
  [\cos^{3}\chi\,\cos{3\theta}] \;
  [\sin^{2}\chi\,\sin{2\varphi}] \nonumber \\
&=& [3(x^2 + y^2) - 4(z^2 + w^2)] \; [x^3 - 3xy^2] \; [2zw].
\label{Rectangular732}
\end{eqnarray}

The fact that each $\widetilde\Psi_{k\ell m}$ may be expressed as a polynomial
proves that the $\widetilde\Psi_{k\ell m}$ are smooth even along the circles
$\chi = 0$ and $\chi = \pi/2$, where the toroidal coordinate system
collapses.

\section{The eigenmodes form a basis}\label{BasisSection}

For each $k$, the set of $\widetilde\Psi_{k\ell m}$ forms a basis for
the space of eigenfunctions on $S^3$ with wave number $k$.
More precisely, define the basis
\begin{equation}
B_k = \{\,\widetilde\Psi_{k\ell m}\quad|\quad|\ell| + |m| \leq k\quad{\rm and}
\quad\ell + m \equiv k\;{\rm (mod\,2)}\,\}.
\label{Eigenbasis}
\end{equation}
To prove that $B_k$ is a basis, we must show that the $\widetilde\Psi_{k\ell m}$
it contains are linearly independent and span the full eigenspace.

{\it Linear independence.}  The inner product of two elements
$\widetilde\Psi_{k\ell m}$ and $\widetilde\Psi_{k\ell'm'}$ of $B_k$ is
\begin{eqnarray}
\label{InnerProduct}
&& \langle \widetilde\Psi_{k\ell m}, \widetilde\Psi_{k\ell'm'} \rangle \nonumber \\
&=& \int_{S^3}\widetilde\Psi_{k\ell m}\widetilde\Psi_{k\ell'm'}\,dV \\
&=& \int_{\chi=0}^{\pi/2} \int_{\theta=0}^{2\pi} \int_{\varphi=0}^{2\pi}
     (\Chi_{k\ell m}\Theta_\ell\Phi_m) \, (\Chi_{k\ell'm'}\Theta_{\ell'}\Phi_{m'})
     \, \cos\chi\,\sin\chi\,d\varphi \, d\theta \, d\chi \nonumber \\
&=& \left(\int_{\chi=0}^{\pi/2}   \Chi_{k\ell m}\Chi_{k\ell'm'}\,
         \cos\chi\,\sin\chi\,d\chi\right)
    \left(\int_{\theta=0}^{2\pi}  \Theta_\ell\Theta_{\ell'}   \, d\theta  \right)
    \left(\int_{\varphi=0}^{2\pi} \Phi_m\Phi_{m'}       \, d\varphi \right). \nonumber
\end{eqnarray}
If $\ell \neq \ell'$ (resp. $m \neq m'$), then the orthogonality of the
circular harmonics $\langle \Theta_\ell,\Theta_{\ell'} \rangle = 0$
(resp. $\langle \Phi_m,\Phi_{m'} \rangle = 0$) immediately implies
$\langle \widetilde\Psi_{k\ell m}, \widetilde\Psi_{k\ell'm'} \rangle = 0$,
proving that $\widetilde\Psi_{k\ell m}$ and $\widetilde\Psi_{k\ell'm'}$ are orthogonal.
Because the $\widetilde\Psi_{k\ell m}$ in $B_k$ are nonzero and pairwise orthogonal,
they must be linearly independent.

{\it Span.}  We have shown that the $\widetilde\Psi_{k\ell m}$ in $B_k$ are
linearly independent.  To prove that they span the full eigenspace,
it suffices to check that the number of elements of $B_k$ equals
the dimension of the full eigenspace, which is known to be $(k+1)^2$.
The set $B_0 = \{\widetilde\Psi_{0,0,0}\}$ has $(0+1)^2 = 1$ element,
and the set $B_1 = \{\widetilde\Psi_{1,+1,0},\widetilde\Psi_{1,-1,0},
\widetilde\Psi_{1,0,+1},\widetilde\Psi_{1,0,-1}\}$ has $(1+1)^2 = 4$ elements,
as required.
For the remaining $B_k$, with $k \geq 2$, we proceed by induction,
assuming that the set $B_{k-2}$ is already known
to contain $((k-2)+1)^2 = (k-1)^2$ elements.
Each element $\widetilde\Psi_{k-2,\ell,m} \in B_{k-2}$ corresponds to
an element $\widetilde\Psi_{k,\ell,m} \in B_{k}$.
The set $B_{k}$ also contains the additional elements
$\widetilde\Psi_{k,0,\pm k}, \widetilde\Psi_{k,\pm 1,\pm (k-1)}, \dots ,
\widetilde\Psi_{k,\pm (k-1),\pm 1}, \widetilde\Psi_{k,\pm k,0}$.
Taking into account the plus-or-minus signs, this gives
$2 + 4 + \dots + 4 + 2 = 2 + 4(k-1) + 2 = 4k$ additional elements.
Adding these to the $(k-1)^2$ elements corresponding to $B_{k-2}$,
we get a total of $(k-1)^2 + 4k = (k+1)^2$ elements, as required.

This completes the proof that $B_k$ is a basis for
the space of eigenfunctions on $S^3$ with wave number $k$.

\section{Normalization}

The $\widetilde\Psi_{k\ell m}$ are already mutually orthogonal
(Section \ref{BasisSection}), so if we normalize them to unit length
they will form an orthonormal basis for the eigenspace.
An orthonormal basis is convenient in cosmological applications,
because it makes it easy to construct an unbiased random
density fluctuation with wave number $k$.

To compute the norm of a given $\widetilde\Psi_{k\ell m}$, set $\ell' = \ell$ and
$m' = m$ in equation (\ref{InnerProduct}), giving
\begin{equation}
\langle \widetilde\Psi_{k\ell m}, \widetilde\Psi_{k\ell m} \rangle
= \left(\int_{\chi=0}^{\pi/2}   \Chi_{k\ell m}^2\,
         \cos\chi\,\sin\chi\,d\chi\right)
    \left(\int_{\theta=0}^{2\pi}  \Theta_\ell^2  \, d\theta  \right)
    \left(\int_{\varphi=0}^{2\pi} \Phi_m^2    \, d\varphi \right).
\label{Norm}
\end{equation}
The $\Theta$ and $\Phi$ integrals are easy to evaluate.
Substituting in the solutions (\ref{ThetaSolution}) and (\ref{PhiSolution})
immediately gives
\begin{equation}
\label{ThetaSolutionGeneric}
   \int_{\theta=0}^{2\pi}  \Theta_\ell^2  \, d\theta  = \pi
 \quad{\rm and}\quad
   \int_{\varphi=0}^{2\pi} \Phi_m^2    \, d\varphi = \pi
\end{equation}
when $\ell$ and $m$ are nonzero, along with the special cases
\begin{equation}
\label{ThetaSolutionSpecial}
   \int_{\theta=0}^{2\pi}  \Theta_0^2  \, d\theta  = 2\pi
 \quad{\rm and}\quad
   \int_{\varphi=0}^{2\pi} \Phi_0^2    \, d\varphi = 2\pi.
\end{equation}
The $\Chi$ integral seems daunting, but luckily the
$\cos^{|\ell|}\chi\,\sin^{|m|}\chi$ factor in the expression
(\ref{RadialSolution}) for $\Chi$ provides exactly the standard
weighting function relative to which the Jacobi polynomials
are normalized!
\begin{eqnarray}
& &  \int_{\chi=0}^{\pi/2} \Chi_{k\ell m}^2\, \cos\chi\,\sin\chi\,d\chi \nonumber \\
&=&  \int_{\chi=0}^{\pi/2} \left[\cos^{|\ell|}\chi\,\sin^{|m|}\chi\;
 P^{(|m|,|\ell|)}_d (\cos 2\chi)\right]^2\, \cos\chi\,\sin\chi\,d\chi \nonumber \\
&=&  \int_{\chi=0}^{\pi/2} (\cos^2\chi)^{|\ell|}\,(\sin^2\chi)^{|m|}\;
\;\left[P^{(|m|,|\ell|)}_d (\cos 2\chi)\right]^2\, \cos\chi\,\sin\chi\,d\chi
\end{eqnarray}
then, changing the variable to $u = \cos 2\chi$,
\begin{equation}
=  \frac{1}{2^{|\ell|+|m|+2}}\;\int_{u=-1}^{+1} (1+u)^{|\ell|}\,(1-u)^{|m|}\;
\;\left[P^{(|m|,|\ell|)}_d (u)\right]^2\, du
\end{equation}
and using the standard normalization for Jacobi polynomials \cite{CRC,AbramowitzStegun},
we get
\begin{equation}
= \frac{(|\ell|+d)!\;(|m|+d)!}{2(k+1)\;d!\;(|\ell|+|m|+d)!}.
\end{equation}

Define the {\it normalized eigenmodes} $\Psi_{k\ell m}$ to be
\begin{equation}
\Psi_{k\ell m}
  = \frac{\widetilde\Psi_{k\ell m}}
  {\sqrt{\langle \widetilde\Psi_{k\ell m}, \widetilde\Psi_{k\ell m} \rangle}}
  = \frac{\widetilde\Psi_{k\ell m}}{\sqrt{2^{\hat{\ell} + \hat{m}}}\;\pi\;\sqrt\frac{(|\ell|+d)!\;(|m|+d)!}{2(k+1)\;d!\;(|\ell|+|m|+d)!}}
\end{equation}
where $\hat\ell$ is 1 when $\ell$ is 0, and $\hat\ell$ is 0 when
$\ell$ is nonzero, and similarly for $\hat m$, to accommodate the
special cases in equations (\ref{ThetaSolutionGeneric}) and
(\ref{ThetaSolutionSpecial}).

The results of this section and the previous one together prove

\vskip0.25cm

\noindent {\bf Theorem 1.}  {\it The $\Psi_{k\ell m}$, taken over
all integers $k$, $\ell$ and $m$ satisfying $|\ell| + |m| \leq k$
and $\ell + m \equiv k$ (mod 2), comprise an orthonormal basis
for the eigenspace of the Laplacian on the 3-sphere.}

\section{The action of an isometry of $S^3$ on the space of eigenmodes}
\label{SectionIsometryAction}

An arbitrary orientation-preserving isometry of the 3-sphere
has matrix
\begin{equation}
\left(
\begin{array}{cccc}
  \cos\Delta\theta & -\sin\Delta\theta &     0       &      0       \\
  \sin\Delta\theta & \quad\cos\Delta\theta &     0       &      0       \\
      0      &      0      & \cos\Delta\varphi & -\sin\Delta\varphi \\
      0      &      0      & \sin\Delta\varphi & \quad\cos\Delta\varphi \\
\end{array}
\right)
\label{IsometryRectangular}
\end{equation}
relative to an appropriate orthonormal $(x,y,z,w)$ coordinate
system on ${\bf R}^4$.  In toroidal coordinates $(\chi,\theta,\phi)$,
the same isometry may be described as
\begin{equation}
\begin{array}{ccl}
\chi    &\to& \chi  \\
\theta  &\to& \theta  + \Delta\theta \\
\varphi &\to& \varphi + \Delta\varphi.
\label{IsometryTorus}
\end{array}
\end{equation}
The action of this isometry on the space of eigenfunctions
is much simpler in toroidal coordinates than in traditional polar coordinates.
For example, if $\ell > 0$ and $-m < 0$, then the isometry
(\ref{IsometryTorus}) maps
\begin{eqnarray}
\Psi_{k,+\ell,-m}(\chi,\theta,\varphi)
 &\to&
    \Psi_{k,+\ell,-m}(\chi,\,\theta  + \Delta\theta,\,\varphi + \Delta\varphi) \nonumber \\
 &=& \Chi_{k\ell m}(\chi)\;\Theta_{+\ell}(\theta  + \Delta\theta)\;\Phi_{-m}(\varphi + \Delta\varphi) \nonumber \\
 &=& \Chi_{k\ell m}(\chi)\;\cos \ell(\theta  + \Delta\theta)\;\sin m(\varphi + \Delta\varphi) \nonumber \\
 &=& \Chi_{k\ell m}(\chi)\;
     (\cos \ell\Delta\theta \cos \ell\theta - \sin \ell\Delta\theta \sin \ell\theta) \nonumber \\
 && \times \quad (\sin m\Delta\varphi \cos m\varphi + \cos m\Delta\varphi \sin m\varphi) \nonumber \\
 &=& \quad \cos \ell\Delta\theta \sin m\Delta\varphi \; \Chi_{k\ell m}(\chi)
          \, \Theta_{+\ell}(\theta) \, \Phi_{+m}(\varphi) \nonumber \\
 &&+ \cos \ell\Delta\theta \cos m\Delta\varphi \; \Chi_{k\ell m}(\chi)
          \, \Theta_{+\ell}(\theta) \, \Phi_{-m}(\varphi) \nonumber \\
 &&- \sin \ell\Delta\theta \sin m\Delta\varphi \; \Chi_{k\ell m}(\chi)
          \, \Theta_{-\ell}(\theta) \, \Phi_{+m}(\varphi) \nonumber \\
 &&- \sin \ell\Delta\theta \cos m\Delta\varphi \; \Chi_{k\ell m}(\chi)
          \, \Theta_{-\ell}(\theta) \, \Phi_{-m}(\varphi) \nonumber \\
 &=& \quad \cos \ell\Delta\theta \sin m\Delta\varphi \;
        \Psi_{k,+\ell,+m}(\chi, \theta, \varphi) \nonumber \\
 &&+ \cos \ell\Delta\theta \cos m\Delta\varphi \;
        \Psi_{k,+\ell,-m}(\chi, \theta, \varphi) \nonumber \\
 &&- \sin \ell\Delta\theta \sin m\Delta\varphi \;
        \Psi_{k,-\ell,+m}(\chi, \theta, \varphi) \nonumber \\
 &&- \sin \ell\Delta\theta \cos m\Delta\varphi \;
        \Psi_{k,-\ell,-m}(\chi, \theta, \varphi) \nonumber
\end{eqnarray}
\begin{equation}
{\rm \dots and\,similarly\,for\,the\,images\,of}\; \Psi_{k,+\ell,+m}, \Psi_{k,-\ell,+m},
  {\rm \,and\,} \Psi_{k,-\ell,-m}.
\label{ActionOnPsi}
\end{equation}
Thus the subspace spanned by the $\Psi_{k,\pm \ell,\pm m}$ is invariant
(setwise but not necessarily pointwise) under the action of the isometry
(\ref{IsometryTorus}).
Typically this subspace is 4-dimensional, but when $\ell$ or $m$
is zero it is 2-dimensional, or only 1-dimensional when both
$\ell$ and $m$ are zero.

Thus the complete eigenspace factors into orthogonal
1-, 2-, and 4-dimensional invariant subspaces, each spanned by
a set $\Psi_{k,\pm \ell,\pm m}$.  To understand the full action of
the isometry, it suffices to understand its action on each
invariant subspace.

\newcommand{\ldt}{\ell\Delta\theta}
\newcommand{\mdp}{m\Delta\varphi}

\begin{description}
  \item[Case 1:]  $\ell = m = 0$.

    The 1-dimensional invariant subspace spanned by $\Psi_{k,0,0}$
    is pointwise fixed by every isometry of the form (\ref{IsometryTorus}).

  \item[Case 2:] $\ell = 0$ or $m = 0$ (but not both).

    For sake of discussion, assume $\ell > 0$ and $m = 0$.
    Relative to the basis
    \begin{eqnarray}
    {\bf d}_1 = \Psi_{k,+\ell,0} \nonumber \\
    {\bf d}_2 = \Psi_{k,-\ell,0}
    \label{NiceBasis2D}
    \end{eqnarray}
    the isometry (\ref{IsometryTorus}) acts as a rotation
    through an angle $\ldt$, with matrix
    \begin{equation}
    \left(
    \begin{array}{rrrr}
       \cos \ldt & \sin \ldt  \\
      -\sin \ldt & \cos \ldt  \\
    \end{array}
    \right).
    \label{Action2D}
    \end{equation}
    A quick computation similar to (\ref{ActionOnPsi}) verifies that
    matrix (\ref{Action2D}) is correct.
    Thus if $\ldt \equiv 0 \;({\rm mod}\, 2\pi)$ the whole subspace
    is fixed pointwise;  otherwise only the origin is fixed.
    Similar conclusions hold when $\ell = 0$ and $m > 0$.

  \item[Case 3:] $\ell > 0$ and $m > 0$.

    Computation (\ref{ActionOnPsi}) expresses the isometry's action on
    each of the $\Psi_{k,\pm \ell,\pm m}$ as a linear combination of the
    $\Psi_{k,\pm \ell,\pm m}$ themselves.  For best results use the rotated basis
    \begin{eqnarray}
    {\bf e}_1 = \sqrt{\frac{1}{2}}\; (\Psi_{k,+\ell,+m} + \Psi_{k,-\ell,-m})
                   = \sqrt{\frac{1}{2}}\; \Chi_{k\ell m}\,\cos(\ell\theta - m\varphi) \nonumber \\
    {\bf e}_2 = \sqrt{\frac{1}{2}}\; (\Psi_{k,-\ell,+m} - \Psi_{k,+\ell,-m})
                   = \sqrt{\frac{1}{2}}\; \Chi_{k\ell m}\,\sin(\ell\theta - m\varphi) \nonumber \\
    {\bf e}_3 = \sqrt{\frac{1}{2}}\; (\Psi_{k,+\ell,+m} - \Psi_{k,-\ell,-m})
                   = \sqrt{\frac{1}{2}}\; \Chi_{k\ell m}\,\cos(\ell\theta + m\varphi) \nonumber \\
    {\bf e}_4 = \sqrt{\frac{1}{2}}\; (\Psi_{k,-\ell,+m} + \Psi_{k,+\ell,-m})
                   = \sqrt{\frac{1}{2}}\; \Chi_{k\ell m}\,\sin(\ell\theta + m\varphi)
    \label{NiceBasis4D}
    \end{eqnarray}
    relative to which the action has matrix
    \begin{equation}
    \left(
    {\small
    \begin{array}{cccc}
     \quad\cos(\ldt - \mdp) &  \sin(\ldt - \mdp) &            0            &          0         \\
       -  \sin(\ldt - \mdp) &  \cos(\ldt - \mdp) &            0            &          0         \\
                0           &          0         &  \quad\cos(\ldt + \mdp) &  \sin(\ldt + \mdp) \\
                0           &          0         &    -  \sin(\ldt + \mdp) &  \cos(\ldt + \mdp) \\
    \end{array}
    }
    \right).
    \label{Action4D}
    \end{equation}
    Clearly the subspace spanned by $\{{\bf e}_1,{\bf e}_2\}$
    (resp. $\{{\bf e}_3,{\bf e}_4\}$) is pointwise
    fixed if and only if $\ldt \equiv \mdp \;({\rm mod}\, 2\pi)$
    (resp. $\ldt \equiv -\mdp \;({\rm mod}\, 2\pi)$).
    If $\ldt \equiv \mdp \equiv 0 \;{\rm or}\; \pi$,
    then the whole 4-dimensional subspace is pointwise fixed.

\end{description}

In summary, an isometry (\ref{IsometryTorus}) fixes
the subspace spanned by
\begin{equation}
\begin{array}{cl}
\{\Psi_{k,0,0}\} & {\rm always} \\
\{\Psi_{k,+\ell,0},\Psi_{k,-\ell,0}\} & {\rm iff}\quad \ldt \equiv 0 \;({\rm mod}\, 2\pi)\\
\{\Psi_{k,0,+m},\Psi_{k,0,-m}\} & {\rm iff}\quad \mdp \equiv 0 \;({\rm mod}\, 2\pi)\\
\begin{array}{l}
\{\sqrt{\frac{1}{2}}\; (\Psi_{k,+\ell,+m} + \Psi_{k,-\ell,-m}), \\
  \quad \sqrt{\frac{1}{2}}\; (\Psi_{k,-\ell,+m} - \Psi_{k,+\ell,-m})\}
\end{array}
   & {\rm iff}\quad \ldt \equiv \mdp \;({\rm mod}\, 2\pi) \\
\begin{array}{l}
\{\sqrt{\frac{1}{2}}\; (\Psi_{k,+\ell,+m} - \Psi_{k,-\ell,-m}), \\
  \quad \sqrt{\frac{1}{2}}\; (\Psi_{k,-\ell,+m} + \Psi_{k,+\ell,-m})\}
\end{array}
   & {\rm iff}\quad \ldt \equiv -\mdp \;({\rm mod}\, 2\pi)
\label{ActionSummary}
\end{array}
\end{equation}
and this is the complete description of the fixed point set.
When working with fixed point sets, please keep in mind that
the eigenmode $\Psi_{k\ell m}$ exists if and only if
$|\ell| + |m| \leq k$ and $\ell + m \equiv k$ (mod 2).

\section{Eigenmodes of lens spaces}
\label{SectionLensSpaces}

The lens space $L(p,q)$ is the quotient of the 3-sphere $S^3$ by
the cyclic group whose generator $g$ is the isometry (\ref{IsometryTorus})
with $\Delta\theta = 2\pi/p$ and $\Delta\varphi = 2\pi q/p$.
Each eigenmode of $L(p,q)$ lifts to a $g$-invariant eigenmode of $S^3$,
and conversely each $g$-invariant eigenmode of $S^3$ projects down
to an eigenmode of $L(p,q)$.  Thus the eigenmodes of $L(p,q)$
correspond to the $g$-invariant eigenmodes of $S^3$.

Substituting $\Delta\theta = 2\pi/p$ and $\Delta\varphi = 2\pi q/p$
into (\ref{ActionSummary}) yields a set of simple integer conditions
showing which eigenmodes of $S^3$ are $g$-invariant:

\vskip0.25cm

\noindent {\bf Theorem 2.}  {\it The eigenspace of the Laplacian
on the lens space $L(p,q)$ has an ortho\-normal basis that, when lifted
to $Z_p$-invariant eigenmodes of the 3-sphere, comprises those
eigenmodes in the left column for which the corresponding
condition in the right column is satisfied, subject to the
restriction that an eigenmode $\Psi_{k\ell m}$ exists
if and only if the integers $k$, $\ell$ and $m$
satisfy $|\ell| + |m| \leq k$ and $\ell + m \equiv k$ (mod 2).}
\begin{equation}
\label{LensSpaceConditions}
\begin{tabular}{|ccc|}
\hline
basis vectors & \quad & condition \\
\hline
$\Psi_{k,0,0}$ & & always \\
$\Psi_{k,+\ell,0},\;\Psi_{k,-\ell,0}$ &  & $\ell  \equiv 0$ (mod $p$)\\
$\Psi_{k,0,+m},\;\Psi_{k,0,-m}$ &  & $qm \equiv 0$ (mod $p$)\\
$\sqrt{\frac{1}{2}}\; (\Psi_{k,+\ell,+m} + \Psi_{k,-\ell,-m}), \quad
  \sqrt{\frac{1}{2}}\; (\Psi_{k,-\ell,+m} - \Psi_{k,+\ell,-m})$
   &  & $\ell \equiv  qm $ (mod $p$) \\
$\sqrt{\frac{1}{2}}\; (\Psi_{k,+\ell,+m} - \Psi_{k,-\ell,-m}), \quad
  \sqrt{\frac{1}{2}}\; (\Psi_{k,-\ell,+m} + \Psi_{k,+\ell,-m})$
&  & $\ell \equiv -qm $ (mod $p$) \\
\hline
\end{tabular}
\end{equation}


\noindent
For a given lens space $L(p,q)$ and wave number $k$,
the number of basis vectors gives the multiplicity of the eigenvalue $k(k+2)$
(see Table \ref{LensSpaceSpectra}).

\begin{table}
\begin{center}
\begin{tabular}{|c|ccccccccccccccc|}
\hline
   $k$   &  0 &  1 &  2 &  3 &  4 &  5 &  6 &  7 &  8 &  9 & 10 & 11 & 12 & 13 & 14 \\
\hline
$S^3$    &  1 &  4 &  9 & 16 & 25 & 36 & 49 & 64 & 81 &100 &121 &144 &169 &196 &225 \\
$L(2,1)$ &  1 &  0 &  9 &  0 & 25 &  0 & 49 &  0 & 81 &  0 &121 &  0 &169 &  0 &225 \\
$L(3,1)$ &  1 &  0 &  3 &  8 &  5 & 12 & 21 & 16 & 27 & 40 & 33 & 48 & 65 & 56 & 75 \\
$L(4,1)$ &  1 &  0 &  3 &  0 & 15 &  0 & 21 &  0 & 45 &  0 & 55 &  0 & 91 &  0 &105 \\
$L(5,1)$ &  1 &  0 &  3 &  0 &  5 & 12 &  7 & 16 &  9 & 20 & 33 & 24 & 39 & 28 & 45 \\
$L(5,2)$ &  1 &  0 &  1 &  4 &  5 &  8 &  9 & 12 & 17 & 20 & 25 & 28 & 33 & 40 & 45 \\
$L(6,1)$ &  1 &  0 &  3 &  0 &  5 &  0 & 21 &  0 & 27 &  0 & 33 &  0 & 65 &  0 & 75 \\
$L(7,1)$ &  1 &  0 &  3 &  0 &  5 &  0 &  7 & 16 &  9 & 20 & 11 & 24 & 13 & 28 & 45 \\
$L(7,2)$ &  1 &  0 &  1 &  2 &  3 &  6 &  7 & 10 & 11 & 14 & 17 & 20 & 25 & 28 & 33 \\
$L(8,1)$ &  1 &  0 &  3 &  0 &  5 &  0 &  7 &  0 & 27 &  0 & 33 &  0 & 39 &  0 & 45 \\
$L(8,3)$ &  1 &  0 &  1 &  0 &  7 &  0 & 11 &  0 & 23 &  0 & 27 &  0 & 45 &  0 & 53 \\
$L(9,1)$ &  1 &  0 &  3 &  0 &  5 &  0 &  7 &  0 &  9 & 20 & 11 & 24 & 13 & 28 & 15 \\
$L(9,2)$ &  1 &  0 &  1 &  2 &  1 &  4 &  7 &  6 &  9 & 14 & 11 & 16 & 21 & 18 & 25 \\
\hline
\end{tabular}
\end{center}
\caption{The number of eigenbasis vectors specified in Theorem 2
  tells the multiplicity of each eigenvalue $k(k+2)$
  in the spectrum of a lens space $L(p,q)$.  Here are the results
  for $p \leq 9$ and $k \leq 14$.}
\label{LensSpaceSpectra}
\end{table}

\section{Eigenmodes of prism spaces}
\label{SectionPrismSpaces}

\begin{figure}
\begin{center}
\epsfig{file=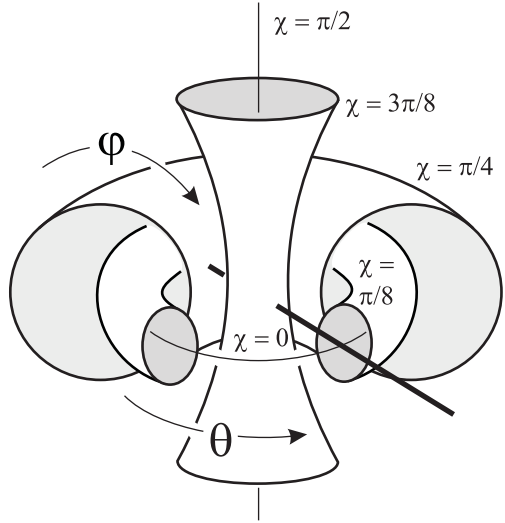}
\caption{The binary dihedral group $D^*_n$ has two generators.
The first generator acts as a lefthanded $\frac{2\pi}{2n}$ corkscrew motion
preserving the toroidal layers just as in the lens space $L(2n,1)$.
The second generator acts as a lefthanded $\frac{2\pi}{4}$ corkscrew motion
along an orthogonal axis (drawn heavy in the figure),
taking, for example, the torus at level $\chi = \frac{\pi}{8}$ to the one
at $\chi = \frac{3\pi}{8}$
and interchanging the roles of $\theta$ and $\varphi$.} \label{fig2}
\end{center}
\end{figure}

The $n^{\rm th}$ {\it prism space} is the quotient $S^3/D^*_n$,
where the binary dihedral group $D^*_n$ is the extension of the
binary cyclic group $Z^*_n = Z_{2n}$ by an order four Clifford
translation along a perpendicular axis (Figure~\ref{fig2}).
The generator of the $Z_{2n}$ action may be written in toroidal
coordinates as
\begin{equation}
   (\chi, \theta, \varphi) \to
      (\, \chi, \quad \theta + \frac{2\pi}{2n}, \quad \varphi + \frac{2\pi}{2n} \,)
\label{IsometryPrismA}
\end{equation}
while the generator of the order four Clifford translation
may be written as
\begin{equation}
   (\chi, \theta, \varphi) \to
      (\, \frac{\pi}{2} - \chi, \quad -\varphi, \quad \pi - \theta \,).
\label{IsometryPrismB}
\end{equation}
The eigenmodes of the prism space naturally correspond to the
$D^*_n$-invariant eigenmodes of $S^3$.  That is, they correspond
to the eigenmodes of $S^3$ that are invariant under both
the $Z_{2n}$ generator (\ref{IsometryPrismA}) and
the $Z_4$ generator (\ref{IsometryPrismB}).
Sections \ref{SectionIsometryAction} and \ref{SectionLensSpaces}
have already determined the action of the $Z_{2n}$ generator
on the space of eigenmodes, with the fixed subspace specified
by Theorem 2 with $p = 2n$ and $q = 1$.
A similar computation shows how the $Z_4$ generator interchanges
the roles of $\ell$ and $m$:
\begin{eqnarray}
\Psi_{k\ell m}(\chi,\theta,\varphi)
 &\to&
     \Psi_{k\ell m}(\pi/2 - \chi,\,-\varphi,\,\pi - \theta) \nonumber \\
 &=& \Chi_{k\ell m}(\pi/2 - \chi)\;\Theta_{\ell}(-\varphi)\;\Phi_{m}(\pi - \theta) \nonumber \\
 &=& \Chi_{k\ell m}(\pi/2 - \chi)\;\Phi_{\ell}(-\varphi)\;\Theta_{m}(\pi - \theta) \nonumber \\
 &=& [\pm\Chi_{km\ell}(\chi)]\;[\pm\Theta_{m}(\theta)]\;[\pm\Phi_{\ell}(\varphi)] \nonumber \\
 &=& \pm \Psi_{km\ell}(\chi,\theta,\varphi).
\label{SecondGeneratorProtoAction}
\end{eqnarray}
The first plus-or-minus sign in the penultimate line of (\ref{SecondGeneratorProtoAction})
will be plus (resp. minus) when $d = \frac{k - (|\ell| + |m|)}{2}$ is
even (resp. odd).  [Proof:  Substituting $\chi \to \frac{\pi}{2} - \chi$ into
(\ref{RadialSolution}) and interchanges the roles of sine and cosine
introduces a factor of $(-1)^d$.]
The third plus-or-minus sign in (\ref{SecondGeneratorProtoAction})
will be plus when $\ell$ is nonnegative, minus otherwise.
The second plus-or-minus sign will be plus when $m$ is either nonnegative and even
or negative and odd, minus otherwise.
Thus we may rewrite
(\ref{SecondGeneratorProtoAction}) as
\begin{eqnarray}
\Psi_{k\ell m}(\chi,\theta,\varphi) &\to&
  [\pm\Chi_{km\ell}(\chi)]\;[\pm\Theta_{m}(\theta)]\;[\pm\Phi_{\ell}(\varphi)] \nonumber \\
 &=& \sigma_{k\ell m} \Psi_{km\ell}(\chi,\theta,\varphi)
\label{SecondGeneratorAction}
\end{eqnarray}
with
$$\sigma_{k\ell m} = (\pm)(\pm)(\pm)(\pm)$$
where
\begin{center}
\begin{tabular}{rcl}
  the first  & $\pm$ & is + if and only if $d$ is even \\
  the second & $\pm$ & is + if and only if $\ell \geq 0$ \\
  the third  & $\pm$ & is + if and only if $m \geq 0$ \\
  the fourth & $\pm$ & is + if and only if $m$ is even.
\end{tabular}
\end{center}

\begin{figure}
\begin{center}
\epsfig{file=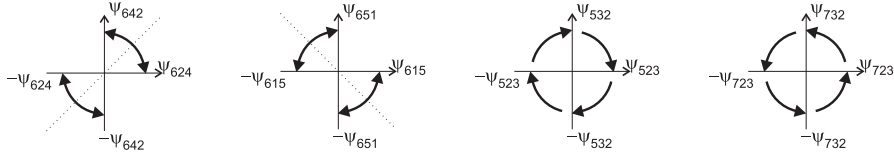,width=12cm}
\caption{
The $Z_4$ generator acts on each 2-dimensional subspace
$\{\Psi_{k\ell m}, \Psi_{km\ell}\}$ as a reflection
(if $\ell$ and $m$ have the same parity) or a rotation
(if $\ell$ and $m$ have opposite parities).} \label{fig3}
\end{center}
\end{figure}

When $\ell \neq m$ the action (\ref{SecondGeneratorAction})
preserves the 2-dimensional plane spanned by
$\Psi_{k\ell m}$ and $\Psi_{km\ell}$.
If the parities of $\ell$ and $m$ agree (both even or both odd)
then the action is either a {\it positive reflection} interchanging
$\Psi_{k\ell m} \leftrightarrow \Psi_{km\ell}$ (Figure~\ref{fig3}a)
or a {\it negative reflection} interchanging
$\Psi_{k\ell m} \leftrightarrow -\Psi_{km\ell}$ (Figure~\ref{fig3}b),
according to the sign of $\sigma_{k\ell m}$.
If the parities of $\ell$ and $m$ disagree (one even and the other odd)
then the action is a {\it quarter turn} taking either
$\Psi_{k\ell m} \rightarrow \Psi_{km\ell} \rightarrow -\Psi_{k\ell m}
\rightarrow -\Psi_{km\ell} \rightarrow \Psi_{k\ell m}$
or the opposite (Figure~\ref{fig3}cd).  The positive and negative
reflections each fix a 1-dimensional line, but the quarter turns
fix only the origin.

When $\ell = m$ the situation is even simpler.
The action (\ref{SecondGeneratorAction}) either fixes
the 1-dimensional line spanned by $\Psi_{k\ell\ell}$
or inverts it, according to whether $\sigma_{k\ell\ell}$ is
positive or negative.

The preceding two paragraphs have found the eigenmodes
fixed by (\ref{SecondGeneratorAction}), i.e. the eigenmodes
invariant under the action of the $Z_4$ generator (\ref{IsometryPrismB}).
To find the eigenmodes of a prism manifold, we must check
which of the modes invariant under the $Z_4$ generator (\ref{IsometryPrismB})
are invariant under the $Z_{2n}$ generator (\ref{IsometryPrismA}) as well.
This is straightforward, because Theorem 2,
with $p = 2n$ and $q = 1$, already tell us
which modes the $Z_{2n}$ generator preserves.
We know {\it a priori} that the eight basis vectors
$\{\Psi_{k,\pm\ell,\pm m},\Psi_{k,\pm m,\pm\ell}\}$
span a space that is setwise invariant under both
the $Z_{2n}$ and the $Z_4$ generators.
Typically this space is 8-dimensional, but its dimension may be less.
Consider the following special cases.  Assume $\ell$ and $m$
are distinct and positive unless otherwise indicated.

\begin{description}
  \item[$\{\Psi_{k00}\}$] (1-dimensional)

      The mode $\Psi_{k00}$ exists if and only if $k$ is even.
      When it exists, it is always fixed by the $Z_{2n}$ generator,
      as indicated in the first line of the conditions
      (\ref{LensSpaceConditions}) in Theorem 2.
      It is fixed by the $Z_4$ generator if and only if
      $\sigma_{k00}$ is positive, which happens if and only if
      $d = k/2$ is even.

  \item[$\{\Psi_{k,\pm\ell,0},\Psi_{k,0,\pm\ell}\}$] (4-dimensional)

      According to the second and third lines of conditions
      (\ref{LensSpaceConditions}), the $Z_{2n}$ generator
      fixes this whole space pointwise when $\ell \equiv 0$ (mod $2n$),
      and otherwise fixes nothing.

      The $Z_4$ generator leaves the two 2-dimensional
      subspaces spanned by $\{\Psi_{k,\ell,0}, \Psi_{k,0,\ell}\}$ and by
      $\{\Psi_{k,-\ell,0},\Psi_{k,0,-\ell}\}$ setwise invariant.
      If $\ell$ is odd, then the $Z_4$ generator acts
      on each 2-dimensional subspace as a quarter turn, fixing
      nothing but the origin.
      If $\ell$ is even (as it must be when $\ell \equiv 0$ (mod $2n$)),
      then the $Z_4$ generator acts on each
      subspace as a positive or negative reflection, fixing
      $\Psi_{k,\ell,0} + \Psi_{k,0,\ell}$ and $\Psi_{k,-\ell,0} - \Psi_{k,0,-\ell}$
      if $\sigma_{k\ell 0}$ is positive (when $d$ is even), or
      $\Psi_{k,\ell,0} - \Psi_{k,0,\ell}$ and $\Psi_{k,-\ell,0} + \Psi_{k,0,-\ell}$
      if $\sigma_{k\ell 0}$ is negative (when $d$ is odd).

  \item[$\{\Psi_{k,\pm\ell,\pm\ell}\}$] (4-dimensional)

      This 4-dimensional subspace factors into two orthogonal
      2-dimensional subspaces, spanned by
      $\{\Psi_{k,\ell,\ell}+\Psi_{k,-\ell,-\ell},
         \Psi_{k,-\ell,\ell}-\Psi_{k, \ell,-\ell}\}$ and
      $\{\Psi_{k,\ell,\ell}-\Psi_{k,-\ell,-\ell},
         \Psi_{k,-\ell,\ell}+\Psi_{k, \ell,-\ell}\}$ respectively,
      each of which is setwise invariant under both the
      the $Z_{2n}$ and the $Z_4$ generators.

      The subspace
      $\{\Psi_{k,\ell,\ell}+\Psi_{k,-\ell,-\ell},
         \Psi_{k,-\ell,\ell}-\Psi_{k, \ell,-\ell}\}$
      is always fixed by the $Z_{2n}$ generator, according
      to the fourth line of conditions (\ref{LensSpaceConditions}).
      It is fixed by the $Z_4$ generator as well if and only if
      $\sigma_{k\ell\ell}$ is positive, which happens if and only if
      $k \equiv 0$ (mod 4).  ($k$ must be even for the mode
      $\Psi_{k\ell\ell}$ to exist, so the only question is whether
      it is equivalent to 0 or 2 modulo 4.)

      The subspace
      $\{\Psi_{k,\ell,\ell}-\Psi_{k,-\ell,-\ell},
         \Psi_{k,-\ell,\ell}+\Psi_{k, \ell,-\ell}\}$
      is fixed by the $Z_{2n}$ generator if and only if
      $2\ell \equiv 0$ (mod $2n$), according
      to the last line of (\ref{LensSpaceConditions}).
      If $\sigma_{k\ell\ell}$ is positive the $Z_4$ generator
      fixes $\Psi_{k,\ell,\ell}-\Psi_{k,-\ell,-\ell}$ but not
      $\Psi_{k,-\ell,\ell}+\Psi_{k, \ell,-\ell}$, while if
      $\sigma_{k\ell\ell}$ is negative the opposite is true.

  \item[$\{\Psi_{k,\pm\ell,\pm m},\Psi_{k,\pm m,\pm\ell}\}$] (8-dimensional)

      This 8-dimensional subspace factors into four orthogonal
      2-dimensional subspaces, spanned respectively by the bases
      \begin{center}
      \begin{tabular}{cl}
          $\{\Psi_{k,+\ell,+m},\Psi_{k,+m,+\ell}\}$ & \\
          $\{\Psi_{k,-\ell,-m},\Psi_{k,-m,-\ell}\}$ & \\
          $\{\Psi_{k,+\ell,-m},\Psi_{k,-m,+\ell}\}$ & \\
          $\{\Psi_{k,-\ell,+m},\Psi_{k,+m,-\ell}\}$ & .
      \end{tabular}
      \end{center}
      If the parities of $\ell$ and $m$ do not match (one even and the other odd)
      then the $Z_4$ generator acts as a quarter turn on each
      of the four subspaces and fixes only the origin.
      If the parities of $\ell$ and $m$ do match (both even or both odd),
      then the $Z_4$ generator acts as a reflection on each of those
      same four subspaces, fixing the following vectors,
      with the choice of signs (consistent as shown) depending on whether
      $\sigma_{k\ell m}$ is positive or negative:
      \begin{center}
      \begin{tabular}{cl}
          $\Psi_{k,+\ell,+m}\pm\Psi_{k,+m,+\ell}$ & \\
          $\Psi_{k,-\ell,-m}\pm\Psi_{k,-m,-\ell}$ & \\
          $\Psi_{k,+\ell,-m}\mp\Psi_{k,-m,+\ell}$ & \\
          $\Psi_{k,-\ell,+m}\mp\Psi_{k,+m,-\ell}$ & .
      \end{tabular}
      \end{center}
      These four vectors comprise a basis for the $Z_4$ generator's
      4-dimensional fixed point space.
      Taking sums and differences of those vectors gives
      a new, more convenient basis for the same space:
      \begin{equation}
      \begin{tabular}{cl}
          $(\Psi_{k,+\ell,+m} + \Psi_{k,-\ell,-m}) \pm (\Psi_{k,+m,+\ell} + \Psi_{k,-m,-\ell})$ & \\
          $(\Psi_{k,+\ell,+m} - \Psi_{k,-\ell,-m}) \pm (\Psi_{k,+m,+\ell} - \Psi_{k,-m,-\ell})$ & \\
          $(\Psi_{k,+\ell,-m} + \Psi_{k,-\ell,+m}) \mp (\Psi_{k,-m,+\ell} + \Psi_{k,+m,-\ell})$ & \\
          $(\Psi_{k,+\ell,-m} - \Psi_{k,-\ell,+m}) \mp (\Psi_{k,-m,+\ell} - \Psi_{k,+m,-\ell})$ & .
      \end{tabular}
      \label{GenericFixedVectors}
      \end{equation}
      The fourth (resp. fifth) line of the conditions (\ref{LensSpaceConditions})
      shows that the $Z_{2n}$ generator fixes the first and fourth
      (resp. second and third) vectors in (\ref{GenericFixedVectors})
      if and only if $\ell \equiv m$ (mod $2n$) (resp. $\ell \equiv -m$ (mod $2n$)).

\end{description}

\noindent
The above cases completely determine the $D^*_n$-invariant
eigenmodes of $S^3$, thus proving

\vskip0.25cm

\noindent {\bf Theorem 3.}  {\it The eigenspace of the Laplacian
on the prism space $S^3/D^*_n$, where $D^*_n$ is the binary dihedral
group of order $4n$, has an orthonormal basis that, when lifted
to $D^*_n$-invariant eigenmodes of the 3-sphere, comprises those
eigenmodes in the left column for which the corresponding
condition in the right column is satisfied, subject to the
restriction that an eigenmode $\Psi_{k\ell m}$ exists
if and only if the integers $k$, $\ell$ and $m$
satisfy $|\ell| + |m| \leq k$ and $\ell + m \equiv k$ (mod 2).}

\begin{equation}
\label{PrismSpaceConditions}
\begin{tabular}{|ccc|}
\hline
basis vectors & \quad & condition \\
\hline
$\Psi_{k,0,0}$ & & $k \equiv 0 \;({\rm mod}\, 4) $\\
\hline
$
\begin{array}{l}
  \sqrt{\frac{1}{2}}\; (\Psi_{k,\ell,0} + \Psi_{k,0,\ell}) \\
  \sqrt{\frac{1}{2}}\; (\Psi_{k,-\ell,0} - \Psi_{k,0,-\ell})
\end{array}
$
& &
$
\begin{array}{c}
\ell \equiv 0\;({\rm mod}\, 2n) \\
\;{\rm and}\; d \;{\rm even}
\end{array}
$\\
\hline
$
\begin{array}{l}
  \sqrt{\frac{1}{2}}\; (\Psi_{k,\ell,0} - \Psi_{k,0,\ell}) \\
  \sqrt{\frac{1}{2}}\; (\Psi_{k,-\ell,0} + \Psi_{k,0,-\ell})
\end{array}
$
 & &
$
\begin{array}{c}
\ell \equiv 0\;({\rm mod}\, 2n) \\
\;{\rm and}\; d \;{\rm odd}
\end{array}
$\\
\hline
$
\begin{array}{l}
  \sqrt{\frac{1}{2}}\; (\Psi_{k,\ell,\ell}+\Psi_{k,-\ell,-\ell}) \\
  \sqrt{\frac{1}{2}}\; (\Psi_{k,-\ell,\ell}-\Psi_{k, \ell,-\ell})
\end{array}
$
 & &
$
\begin{array}{c}
k \equiv 0 \;({\rm mod}\, 4)
\end{array}
$\\
\hline
$
\sqrt{\frac{1}{2}}\; (\Psi_{k,\ell,\ell}-\Psi_{k,-\ell,-\ell})
$
  & &
$
\begin{array}{c}
2\ell \equiv 0\;({\rm mod}\, 2n) \\
\;{\rm and}\; k \equiv 0\;({\rm mod}\, 4)
\end{array}
$\\
\hline
$
\sqrt{\frac{1}{2}}\; (\Psi_{k,-\ell,\ell}+\Psi_{k, \ell,-\ell})
$
 & &
$
\begin{array}{c}
2\ell \equiv 0\;({\rm mod}\, 2n) \\
\;{\rm and}\; k \equiv 2\;({\rm mod}\, 4)
\end{array}
$\\
\hline
$
\begin{array}{l}
  \frac{1}{2}\; ((\Psi_{k,+\ell,+m} + \Psi_{k,-\ell,-m}) + (\Psi_{k,+m,+\ell} + \Psi_{k,-m,-\ell})) \\
  \frac{1}{2}\; ((\Psi_{k,+\ell,-m} - \Psi_{k,-\ell,+m}) - (\Psi_{k,-m,+\ell} - \Psi_{k,+m,-\ell}))
\end{array}
$
 & &
$
\begin{array}{c}
\ell \equiv m\;({\rm mod}\, 2n) \\
     \;{\rm and}\; \sigma_{k\ell m} > 0
\end{array}
$\\
\hline
$
\begin{array}{l}
  \frac{1}{2}\; ((\Psi_{k,+\ell,+m} + \Psi_{k,-\ell,-m}) - (\Psi_{k,+m,+\ell} + \Psi_{k,-m,-\ell})) \\
  \frac{1}{2}\; ((\Psi_{k,+\ell,-m} - \Psi_{k,-\ell,+m}) + (\Psi_{k,-m,+\ell} - \Psi_{k,+m,-\ell}))
\end{array}
$
 & &
$
\begin{array}{c}
\ell \equiv m\;({\rm mod}\, 2n) \\
     \;{\rm and}\; \sigma_{k\ell m} < 0
\end{array}
$\\
\hline
$
\begin{array}{l}
  \frac{1}{2}\; ((\Psi_{k,+\ell,+m} - \Psi_{k,-\ell,-m}) + (\Psi_{k,+m,+\ell} - \Psi_{k,-m,-\ell})) \\
  \frac{1}{2}\; ((\Psi_{k,+\ell,-m} + \Psi_{k,-\ell,+m}) - (\Psi_{k,-m,+\ell} + \Psi_{k,+m,-\ell}))
\end{array}
$
 & &
$
\begin{array}{c}
\ell \equiv -m\;({\rm mod}\, 2n) \\
     \;{\rm and}\; \sigma_{k\ell m} > 0
\end{array}
$\\
\hline
$
\begin{array}{l}
  \frac{1}{2}\; ((\Psi_{k,+\ell,+m} - \Psi_{k,-\ell,-m}) - (\Psi_{k,+m,+\ell} - \Psi_{k,-m,-\ell})) \\
  \frac{1}{2}\; ((\Psi_{k,+\ell,-m} + \Psi_{k,-\ell,+m}) + (\Psi_{k,-m,+\ell} + \Psi_{k,+m,-\ell}))
\end{array}
$
 & &
$
\begin{array}{c}
\ell \equiv -m\;({\rm mod}\, 2n) \\
     \;{\rm and}\; \sigma_{k\ell m} < 0
\end{array}
$\\
\hline
\end{tabular}
\end{equation}

Note that when $k$ is odd, none of the above conditions are satisfied,
and the eigenbasis is empty.  This is not surprising, because
when $k$ is odd the eigenmodes correspond to odd-degree homogeneous
polynomials (Section \ref{SectionTrigAndPolynomialForms}),
which are anti-symmetric under the action of the antipodal map,
and all groups $D^*_n$ contain the antipodal map.

When $k$ is even, the total number of eigenmodes for given
$D^*_n$ and $k$ agrees with the multiplicities given by Ikeda's formulas
$(2\bar k+1)([\bar k/n] + 1)$ (for $\bar k$ even) and
$(2\bar k+1)[\bar k/n]$ (for $\bar k$ odd) from Theorem 4.3 of \cite{Ikeda},
where $\bar k = k/2$ is half the wave number and $[\bar k/n]$ denotes
the integer part of $\bar k/n$.

\section{Acknowledgments}

We thank Dorian Goldfeld for his generous help evaluating our sums.
JRW thanks the MacArthur Foundation for its support.


\begin{thebibliography}{99}

\bibitem{LuminetStarkmanWeeks}
J.-P. Luminet, G.D. Starkman and J.R. Weeks, ``Is space finite?'',
Scientific American, April 1999 [{\tt
http://www.sciam.com/1999/0499issue/0499weeks.html}].

\bibitem{CornishWeeks}
N.J. Cornish and J.R. Weeks, ``Measuring the shape of the
universe'', Notices of the American Mathematical Society {\bf 45}
(1998) 1463-1471 [{\tt astro-ph/9807311}].

\bibitem{Luminet}
J.-P. Luminet, {\it L'Univers chiffonn\'e}, Fayard, Paris, 2001.

\bibitem{CornishSpergelStarkman}
N.J. Cornish, D.N. Spergel and G.D. Starkman, ``Circles in the
sky: finding topology with the microwave background radiation'',
Class. Quantum Grav. {\bf 15} (1998) 2657-2670 [{\tt
astro-ph/9801212}].

\bibitem{LLL}
R. Lehoucq, M. Lachi\`eze-Rey and J.-P. Luminet,
``Cosmic crystallography'',
Astron. Astrophys. {\bf 313} (1996) 339-346
[{\tt gr-qc/9604050}].

\bibitem{JaffeEtAl}
A.H. Jaffe {\it et al.}, ``Cosmology from MAXIMA-1, BOOMERANG and
COBE/DMR cosmic microwave background observations'', Phys. Rev.
Lett. {\bf 86} (2001) 3475-3479 [{\tt astro-ph/0007333}].

\bibitem{GLLUW}
E. Gausmann, R. Lehoucq, J.-P. Luminet, J.-P. Uzan and J. Weeks,
``Topological lensing in spherical spaces'', Class. Quantum Grav.
{\bf 18} (2001) 5155-5186 [{\tt gr-qc/0106033}].

\bibitem{Levin}
J. Levin, ``Topology and the cosmic microwave background'', Phys.
Rep. {\bf 365} (2002) 251-333 [{\tt gr-qc/0108043}].

\bibitem{IkedaYamamoto}
A. Ikeda and Y. Yamamoto,
``On the spectra of 3-dimensional lens spaces'',
Osaka J. Math. {\bf 16} (1979) 447-469.

\bibitem{BergerGauduchonMazet}
M. Berger, P. Gauduchon and E. Mazet,
{\it Le spectre d'une vari\'et\'e riemannienne},
Lecture Notes in Math. {\bf 194}, Springer-Verlag (1971).

\bibitem{SteinWeiss}
E.M. Stein and G. Weiss, {\it Introduction to Fourier analysis on
Euclidean spaces}, Princeton University Press (1971).

\bibitem{CRC}
S.M. Selby (Ed.), {\it CRC Standard Mathematical Tables}, $19^{\rm
th}$ edition, Chemical Rubber Company, Cleveland, 1971.

\bibitem{AbramowitzStegun}
M. Abramowitz and I.A. Stegun (Eds.), {\it Handbook of
Mathematical Functions with Formulas, Graphs, and Mathematical
Tables}, Dover Publications (1972).

\bibitem{Ikeda}
A. Ikeda, ``On the spectrum of homogeneous spherical space
forms'', Kodai Math. J. {\bf 18} (1995) 57-67.

\end{thebibliography}
\end{document}